\documentclass[11pt]{amsart}

\title{Undecidability of the submonoid membership problem for a  sufficiently large finite direct power of the Heisenberg group}

\author{Vitaly Roman'kov}
\address{Sobolev Institute of Mathematics, SB RAS, Russia}
\curraddr{}
\email{romankov48@mail.ru}

\usepackage{amsmath,amsthm,amsfonts,amssymb}
\usepackage{graphicx}
\usepackage{hyperref}
\usepackage[usenames]{color}

\newtheorem{theorem}{Theorem}[section]
\newtheorem{lemma}[theorem]{Lemma}
\theoremstyle{definition}

\newtheorem{proposition}[theorem]{Proposition}

\newcounter{comcount}

\date{}

\begin{document}

\maketitle

\begin{abstract}
The submonoid membership problem for a finitely generated group $G$ is the decision problem, where for a given finitely generated submonoid $M$ of $G$ and a group element  $g$ it is asked whether $g \in M$. In this paper, we prove that for a sufficiently large direct power $\mathbb{H}^n$
 of the Heisenberg group $\mathbb{H}$, there exists a finitely generated submonoid $M$ whose membership problem is algorithmically unsolvable. Thus, an answer is given to the question of M. Lohrey and B. Steinberg about the existence of a finitely generated nilpotent group with an unsolvable submonoid membership problem. It also answers the question of T. Colcombet, J. Ouaknine, P. Semukhin and J. Worrell about the existence of such a group in the class of direct powers of the Heisenberg group. This result implies the existence of a similar submonoid in any free nilpotent group $N_{k,c}$ of sufficiently large rank $k$ of the class $c\geq 2$. The proofs are based on the undecidability of Hilbert's 10th problem and interpretation of Diophantine equations in nilpotent groups.
 

\end{abstract}

\maketitle

 \section{Introduction}
\label{sec:1}
 
The submonoid membership problem for finitely generated nilpotent groups, which has attracted the attention of a number of researchers in recent years, is considered.  Recall that this is the problem of the existence of an algorithm that determines, given an arbitrary element $g$ and a finitely generated submonoid $M$ of a group $G$, whether $g$ belongs to $M$. Note that in \cite{RomTwo} the author announced a negative solution to this problem for a free nilpotent group $N_{k,c}$ of nilpotency class $c$ at least two of sufficiently large rank $k$. This result will appear in \cite{RomIzv}. This gives an answer to the well-known question of M. Lohrey and B. Steinberg (\cite{LohSur}, Open problem 24)  about the existence of a finitely generated nilpotent group with an unsolvable  submonoid membership problem. Moreover, the existence in $N_{k,c}$ of a finitely generated submonoid $M$  with the unsolvable  membership problem was established. The proof shows how, from an arbitrary Diophantine equation $P$, an element $g$ and a finitely generated submonoid $M$ of the group $N_{k,c}$ are effectively constructed such that $g$ belongs to $M$ if and only if the equation $P$ is solvable in integers. Then the undecidability of Hilbert's 10th problem allows us to obtain from this result the undecidability of the membership problem for $N$ with respect to $M$.

In \cite{Col} the authors prove that the somewhat more general  the subsemigroup membership problem is decidable for the Heisenberg group $\mathbb{H}=H(3, \mathbb{Z})$ consisting of upper triangular integer matrices with units along the diagonal. In other words, the Heisenberg group is a free nilpotent group of rank two and class two. Earlier in \cite{Ko} it was shown how to solve the  problem of belonging of the identity matrix to finitely generated subsemigroups in $\mathbb{H}$.
 In \cite{Col}, the question was raised about the decidability of the submonoid membership problem for a direct power of the Heisenberg group.
 
This paper is a continuation of our previous papers: the paper \cite{RomTwo} mentioned above and the recently published paper \cite{RomPos}, in which sufficient conditions for the solvability of the submonoid  membership problem for a free nilpotent group of the  class $2$ with respect to a given submonoid $M$ were presented, as well as the upcoming paper \cite{RomIzv}. Our objective in this paper is to prove the undecidability of the submonoid membership problem for a sufficiently large finite direct power $\widetilde{\mathbb{H}}=\mathbb{H}^n$ of the Heisenberg group $\mathbb{H}$. Just as in \cite{RomIzv}, we prove that, given any Diophantine equation $P$, one can effectively construct an element $g$ and a finitely generated submonoid $M$ of  the group $\widetilde{\mathbb{H}}$ such that  $g$ belongs to $M$ if and only if the equation $P$ is solvable in integers. Again, the undecidability of Hilbert's 10th problem allows us to obtain from this result the undecidability of the membership problem for $\widetilde{\mathbb{H}}$ with respect to $M$.

This result also easily implies the existence of a submonoid $\widetilde{M}$ of a free nilpotent group $N_{k,c}$ of sufficiently large rank $k$ of the class $c \geq 2$ with the unsolvable problem of belonging to the submonoid $\widetilde{M }$.

Note that a special case of the submonoid membership problem is the classical membership problem, which came from M. Dehn, where $M$ is a finitely generated subgroup. In a different terminology, it is called the generalized word problem.  A.I. Maltsev \cite{Mal} showed that this problem is decidable for any finitely generated nilpotent group. It is worth noting that in the class of finitely generated nilpotent groups, almost all basic algorithmic problems (word, conjugacy, isomorphism, etc.) are solved positively. See surveys \cite{NRR}, \cite{RemRom}, \cite{RomAT}. The exceptions are the problem of endomorphic reducibility and a number of problems related to equations and identities in groups. See \cite{RomEnd}, \cite{RomEq}, \cite{RomDio}, \cite{Kle}, \cite{Rep} on this subject.

The submonoid membership problem is the most important fragment of the more general  rational subset membership problem. See survey \cite{LohSur}.

The submonoid membership problem for a non-commutative group is currently considered as a transfer of the classical problem of integer linear programming, where the submonoid membership problem   for a free abelian group appears on a non-commutative platform. A new line of research has emerged and is being developed -- noncommutative discrete optimization. The chapter ``Discrete optimization in groups'' in the book \cite{GAGTA} is devoted to this direction. In this case, special attention is paid to the class of finitely generated nilpotent groups, which is closest to the class of abelian groups.

\section{Diophantine equations and Skolem systems}
\label{sec:2}

Let $\zeta_1, \ldots , \zeta_t$ be an arbitrary set of commuting variables. A polynomial $D(\zeta_1, \ldots , \zeta_t)$ with integer coefficients in these variables is called {\it Diophantine}.

In this paper, we will write an arbitrary Diophantine equation in the form
\begin{equation}
\label{eq:1}
D(\zeta_1, \ldots , \zeta_t) = \upsilon ,\, \upsilon \in \mathbb{Z},
\end{equation} 
\noindent
where the polynomial from the left side has zero constant term.

\subsection
 {Skolem systems}
 \label{ssec:1}
 In the monograph \cite{Skolem}, T. Skolem showed that any Diophantine equation  is equivalent to a system of equations in a larger number of variables of three types:
$\zeta \zeta' - \zeta'' =0,\, \zeta + \zeta' - \zeta''=0$\, and $\zeta - \zeta' = 0 ,$\, as well as one equation of the form $\zeta - \upsilon =0\, (\upsilon \in \mathbb{Z})$, where each variable  $\zeta , \zeta', \zeta'' $ either occurs in the notation of the original polynomial $D(\zeta_1, \ldots , \zeta_t)$, or is introduced additionally. Such a system is called the {\it Skolem system}. In what follows, we also write the equations of the Skolem system in the form $\zeta\zeta' = \zeta''$,\,
$\zeta + \zeta' = \zeta''$,\, $\zeta = \zeta',$\, and\, $\zeta = \upsilon ,$\, respectively.

 \subsubsection
 {Algorithm for obtaining the Skolem system}
 \label{sssec:1} 
 We will show how to write the Skolem system using an equation of the form (\ref{eq:1}). We assume that either all the coefficients of the polynomial $D(\zeta_1, \ldots , \zeta_{t})$ are positive, or there are coefficients of different signs among them. If initially all these coefficients are negative, then we pass to the equation obtained by multiplying both parts of (\ref{eq:1}) by $-1$.
 
 Let us take one of the non-linear monomials on the left side of the considered equation. Let $\zeta \zeta'$ be the product of its two factors. Let us introduce a new variable $\zeta''$ and write a new equation $\zeta\zeta'= \zeta''$ into the system, simultaneously replacing the product $\zeta \zeta'$ in the monomials  by $\zeta''$. The degree of the monomial under consideration will decrease by one. If it has become linear, go to the next monomial. If not, then we continue to act similarly with the given monomial until it becomes linear.
 
 Then we move on to the next monomial, and so on. As a result, the left side of the equation will be represented as an algebraic sum of variables. Next, we introduce new variables, replacing $\zeta + \zeta'$ in this sum by $\zeta''$ (similarly, $- \zeta - \zeta'$ is replaced by $- \zeta''$ ), adding the equation $\zeta + \zeta' = \zeta''$ to the system. We continue this process. If all coefficients in the algebraic sum are equal to $1$, then the last equation will be of the form $\zeta = \upsilon$, which we will also include in the system. If terms of different signs were present, then by transforming all the terms on the left side, we arrive at an equation of the form $\zeta - \zeta' = \upsilon$. Then we set $\zeta = \zeta' + \zeta''$ and $\zeta'' = \upsilon .$

 \subsubsection{Nonnegative  Diophantine equations and Skolem systems}
 \label{sssec:2}
 A Diophantine equation that is considered solvable if it has a solution in nonnegative integers is called {\it nonnegative}. Similarly, a Skolem system for which decidability means the existence of a solution in nonnegative integers is called {\it nonnegative}.

 \begin{lemma} 
 \label{le:1}
  Solvability of an arbitrary Diophantine equation (\ref{eq:1})
is equivalent to the solvability of some nonnegative Diophantine equation in $2t$ variables effectively constructed from this equation. The resulting equation is equivalent to the nonnegative Skolem system $S_{\upsilon}$.
\end{lemma}
Proof. 
 We write each variable $\zeta_i$ as the difference of the new variables $\zeta_i' - \zeta_i''.$ Substituting these differences for the variables of the equation (\ref{eq:1}), we obtain the nonnegative Diophantine equation 
\begin{equation}
\label{eq:2}
D_1(\zeta_1', \zeta_1'', \ldots , \zeta_t', \zeta_t'') = 0.
\end{equation}
Obviously, the solvability of the equation (\ref{eq:1}) in integers implies the solvability of the equation (\ref{eq:2}) in nonnegative integers, and vice versa.

Based on the nonnegative equation obtained in this way, we build the Skolem system, as described in the  subsubsection \ref{sssec:1}. All substitutions of the form $\zeta \zeta ' = \zeta ''$ and $\zeta + \zeta ' = \zeta ''$ lead to nonnegative variables of the Skolem system.  
An exception is possible only at the final replacement, when it is necessary to transform the equation of the form $\zeta - \zeta' = \upsilon$   for $\upsilon < 0$. Then we set $\zeta ' = \zeta + \zeta ''$ and $\zeta'' = - \upsilon .$ In all cases the last equation has the form $\zeta = |\upsilon |.$
\hfill $\square$

Consider the obtained nonnegative Skolem system $S_{\upsilon}$.  For what follows, we need the renumbering of variables, the introduction of new variables, and the ordering of the equations of the $S_{\upsilon}$ system. For simplicity, the notation $S_{\upsilon}$ does not change in this process.

Assume that $S_{\upsilon}$ contains $e$ equations of the form $\zeta_i\zeta_j=\zeta_l$. Introducing new variables $\zeta'$ and making appropriate substitutions of the form $\zeta_i$ for $\zeta'$, we achieve that each variable will appear in these equations exactly once. Equations of the form $\zeta'=\zeta_i$ will be added to the $S_{\upsilon}$ system. Next, we renumber the variables in such a way that all $e$ equations of the indicated form take the form
$$
\zeta_{1}\zeta_{2}= \zeta_{3},
$$
 \begin{equation}
\label{eq:3} 
  \ldots , \\ 
\end{equation}
$$
 \zeta_{3(e-1)+1} \zeta_{3(e-1)+2}= \zeta_{3e}.
 $$
Let the system $S_{\upsilon}$ contains $d$ equations of the form $\zeta_i+\zeta_j = \zeta_l.$ Similarly to the case just considered, we will ensure that among the variables of the considered set of equations there will be no variables of the previous subsystem, and each variable in their entries will appear in these equations exactly once. Next, we renumber the variables of this subsystem in such a way that all $d$ equations of the indicated form will include only the variables $\zeta_{3e+1}, \ldots , \zeta_{3(e+d)},$ and the subsystem itself will take the form
$$
\zeta_{3e+1} + \zeta_{3e+2} = \zeta_{3(e+1)},
$$
\begin{equation}
\label{eq:4}
 \ldots ,
 \end{equation} 
 $$\zeta_{3(e+d-1) +1} + \zeta_{3(e+d-1) +2} = \zeta_{3(e+d)}.$$

Next, we write the third system, consisting of equations related to the equality of variables.  Let us write  all equalities of the form $\zeta_i=\zeta_j,$ for pairs with different indices $i, j\leq e+d$, which follow from the set of all equalities. Moreover, it suffices to write down a subsystem in which each variable occurs exactly once. Let's renumber all the equations of this subsystem    by assigning them the numbers $e+d +1, \ldots , e + d + q$, respectively. We have the system of equations 
\begin{equation}
\label{eq:5}
P_k \sim \zeta_{i(k)} = \zeta_{j(k)},\, i(k)\not= j(k), i(k), j(k)\leq 3(e+d), k= e+d+1, \ldots , e+d+q. 
\end{equation}
It remains to write a special equation
\begin{equation}
\label{eq:6}
\zeta_t = |\upsilon |.
\end{equation}

Except for the trivial equation (\ref{eq:1}) of the form $\zeta_1=\upsilon \, (t=1)$, the variable $\zeta_t$ is present in the system (\ref{eq:5}). Hence, $\zeta_1, \ldots , \zeta_{3e}, \zeta_{3e+1}, \ldots , \zeta_{3(e+d)}$ are all variables of the system $S_{\upsilon}$. 
 It is obvious that the system $S_{\upsilon }$ is equivalent to the new system  thus replaced.

\section{Auxiliary assertions}
\label{sec:3}

Now we prove a number of auxiliary assertions. The commutator $[g, f]$ of two elements $g$ and $f$ of a group $G$ is defined as $g^{-1}f^{-1}gf$. Then $gf = fg[g,f]$.

Recall that the group $\mathbb{H}$ is generated by the transvections $a=t_{12}$ and $b=t_{23}$, and its center is the infinite cyclic group generated by their commutator $c=[b, a]=t_{13}^{-1}$. Then for any $\alpha , \beta \in \mathbb{Z},$  $b^{\beta}a^{\alpha} =a^{\alpha}b^{\beta}c^{\alpha \beta}.$

\begin{lemma}
\label{le:2}
Let $M$ be a submonoid of $\mathbb{H}$ generated by $g_1=ac, b$ and $g_2=a^{-1}$. 
Then  any representation of $b$ in terms of the generators of $M$ has the form 
\begin{equation}
\label{eq:7}
b = g_1^{\zeta}bg_2^{\zeta}, \zeta \in \mathbb{N}\cup\{0\}.
\end{equation}
\end{lemma}
Proof. 
Obviously, the generator $b$ of the submonoid $M$ appears exactly once among the factors of the right-hand side (\ref{eq:7}), and the total exponents of the occurrences of the other two generators $g_1$ and $g_2$ are equal. Then the cancellation of the degree $c^{\zeta}$ occurs only for the indicated arrangement of the generators of the submonoid on the right-hand side (\ref{eq:7}).

The cancellation  occurs during the commutator collecting process, which consists in the transition of $g_2^{\zeta}$ through $b$. Namely, 
\begin{equation}
\label{eq:8}
g_1^{\zeta}bg_2^{\zeta} = (a^{\zeta}c^{\zeta})a^{-\zeta}b[b, a^{-\zeta}]  = b.
\end{equation}
 \hfill $\square$
 
 The scheme of the exact location of the generators of the submonoid $M$ when expressing the element $b$ is as follows.
 $$
 \begin{vmatrix}
 &&g_1^{\zeta}& b & g_2^{\zeta}\\
 \mathbb{H}:&b=&(ac)^{\zeta}& b & a^{-\zeta}
 \end{vmatrix}.
 $$
  \begin{center} Table 1.
 \end{center}
  
In the following lemmas, $\mathbb{H}^k$ denotes the direct product of $k$ copies of the group $\mathbb{H}$. Denote the transvections $t_{12}, t_{23}, t_{13}^{-1}$  in the $i$-th copy ($i=1, \ldots , k $) as $a_i, b_i, c_i$ respectively. 

The following lemma allows us to interpret equations of the form $\zeta + \zeta ' = \zeta ''$ in the group $\mathbb{H}^4$.
\begin{lemma} 
\label{le:3}
  Let $M$ be a submonoid of $\mathbb{H}^4$ generated by $g_1=a_1c_1c_4, g_2=a_2c_2c_4, g_3=a_3c_3c_4^{-1}, g_4=a_1^{-1}, g_5=a_2^{-1}, g_6=a_3^{-1}$, and $f_1 =b_1b_2b_3.$ Then the representation of 
$b_{1-3}$ in terms of the generators of $M$ has the form
\begin{equation}
\label{eq:9}b_{1\to 3}=g_1^{\zeta}g_2^{\zeta '}g_3^{\zeta ''}f_1g_4^{\zeta}g_5^{\zeta '}g_6^{\zeta ''}\end{equation}
\noindent provided that $\zeta + \zeta ' = \zeta ''.$ For given positive $\zeta , \zeta ', \zeta ''$, the form (\ref{eq:9}) is uniquely determined up to a permutation of the factors $g_i$\, $(i=1, 2, 3)$ on the left side and $g_j $\, $(j = 4 ,5,6)$ on the right side of the factor $f_1$. For null value of $\zeta , \zeta '$ or $\zeta ''$, you can also assume that the corresponding generator  is located as indicated.
\end{lemma}
Proof. 
Obviously, the generator $f_1$ of the submonoid $M$  occurs exactly once  among the factors of the right-hand side of (\ref{eq:9}),  and the total exponents of  occurrences of any pair of generators $g_i$ and $g_{i+3}$ for $i = 1, 2, 3$   are the same, say $\zeta , \zeta '$ and $\zeta ''$, respectively. The location of the factors $g_i$ and $g_{i+3}$  for $i=1,2,3$ with respect to $f_1$ follows from Lemma \ref{le:2}. Then the equality  
$\zeta + \zeta ' = \zeta ''$ is necessary and sufficient for the indicated occurrence of the element $b$ in the submonoid $M$.
 \hfill $\square$
 
Let $\mathbb{H}^4 = \mathbb{H}(1) \times \ldots \times \mathbb{H}(4).$ The scheme of the exact location of the components of the generators of the submonoid $M$ when expressing the element $b_{1\to 3}$ is as follows (empty positions correspond to trivial elements).
 $$
 \begin{vmatrix}
 &b_{1\to 3}=&g_1^{\zeta}&g_2^{\zeta '}&g_3^{\zeta ''}&f_1&g_4^{\zeta}&g_5^{\zeta '}&g_6^{\zeta ''}\\
\mathbb{H}(1):&b_1=&a_1^{\zeta}c_1^{\zeta '}&&&b_1&a_1^{-\zeta}&\\
\mathbb{H}(2):&b_2=&&a_2^{\zeta '}c_2^{\zeta '}&&b_2&&a_2^{-\zeta '}\\
\mathbb{H}(3):&b_3=&&&a_3^{\zeta ''}c_3^{\zeta ''}&b_3&&&a_3^{-\zeta ''}\\
\mathbb{H}(4):&1=&c_4^{\zeta}&c_4^{\zeta '}&c_4^{-\zeta ''}&&&\\
 \end{vmatrix}
 $$
 \begin{center} Table 2.
 \end{center}
  
\begin{lemma}
\label{le:4}
Let $M$ be a submonoid of $\mathbb{H}^6$ generated by $g_1 = a_1c_1, g_2 = a_2c_2, g_3 = a_1^{-1}a_3c_3, g_4= a_2^{-1}a_4c_4, f_1 = b_1b_2, f_2 = b_3b_4,g_5= a_3^{-1}a_5c_5, 
g_6=a_4^{-1}a_6c_6, f_3 = b_5b_6, g_7 = a_5^{-1}, g_8 = a_6^{-1}.$   Then the representation of 
$b_{1-6}= b_1 \cdot \ldots \cdot b_6$  in terms of the generators of $M$ has the form
\begin{equation}
\label{eq:10}
b_{1-6}=g_1^{\zeta}g_2^{\zeta '}f_1g_3^{\zeta}g_4^{\zeta '}f_2g_5^{\zeta}g_6^{\zeta '}f_3g_7^{\zeta}g_8^{\zeta '}.\end{equation}
 For given positive $\zeta , \zeta '$ the form (\ref{eq:10}) is defined uniquely up to a permutation of the  generators  $g_1, g_2$ on the left side  and $g_3,g_4$ on the right side of $f_1$, $g_3,g_4$ on the left side and $g_5,g_6$ on the right side of $f_2$, $g_5,g_6$ on the the left side  and $g_7,g_8$ on the right side of $f_3$. For null value of $\zeta $ or  $\zeta '$, you can also assume that the corresponding generator  is located as indicated.

\end{lemma}
Proof.
The proof is similar to the proof of Lemma \ref{le:3} and follows from Lemma \ref{le:2}. Obviously, each of the generators $f_1, f_2, f_3$ of the submonoid $M$  occurs exactly once  among the factors of the right-hand side of (\ref{eq:9}),  and the total exponents of  occurrences of any pair of generators $(g_1, g_3), (g_2, g_4), (g_3,g_5), (g_4,g_6)$ and $(g_5,g_7), (g_6,g_8)$    are the same, respectively.  Then the quadruples of generators $(g_1, g_3, g_5, g_7)$ and $g_2, g_4, g_6, g_8)$ have the same degrees, say $\zeta$ and  $\zeta '$, respectively. The location of the factors $g_i$th  with respect to $f_j$th  follows from Lemma \ref{le:2}. 
 \hfill $\square$

Let $\mathbb{H}^6 = \mathbb{H}(1) \times \ldots \times \mathbb{H}(6).$ The scheme of the exact location of the components of the generators of the submonoid $M$ when expressing the element $b_{1\to 6}$ is as follows (empty positions correspond to trivial elements).
 $$
 \begin{vmatrix}
 &b_{1\to 6}=&g_1^{\zeta}g_2^{\zeta '}&f_1&g_3^{\zeta}g_4^{\zeta '}&f_2&g_5^{\zeta}g_6^{\zeta '}&f_3&g_7^{\zeta}g_8^{\zeta '}\\
\mathbb{H}(1):&b_1=&a_1^{\zeta}c_1^{\zeta}&b_1&a_1^{-\zeta}   &&&&\\
\mathbb{H}(2):&b_2=&a_2^{\zeta '}c_2^{\zeta '}&b_2&a_2^{-\zeta '}&&&&\\
\mathbb{H}(3):&b_3=&&&a_3^{\zeta}c_3^{\zeta}&b_3&a_3^{-\zeta}&&\\
\mathbb{H}(4):&b_4=&&&a_4^{\zeta '}c_4^{\zeta '}&b_4&a_4^{-\zeta '}&&\\  
\mathbb{H}(5):&b_5=&&&&&a_5^{\zeta}c_5^{\zeta}&b_5&a_5^{-\zeta}\\  
\mathbb{H}(6):&b_6=&&&&&a_6^{\zeta '}c_6^{\zeta '}&b_6&a_6^{-\zeta '}\\  
 \end{vmatrix}.
 $$
 \begin{center} Table 3.
 \end{center}
 
 \begin{lemma}
\label{le:5}
Consider a group $\mathbb{H}^8$ in which the first $6$  components form the group  $\mathbb{H}^6$ from Lemma \ref{le:4}. Let $M$ be a submonoid of $\mathbb{H}^8$ generated by $g_1 = a_1c_1a_7, g_2=a_2c_2, g_3= a_1^{-1}a_3c_3, g_4= a_2^{-1}a_4c_4b_7, f_1=b_1b_2, f_2=b_3b_4, g_5= a_3^{-1}a_5c_5a_7^{-1}, 
g_6=a_4^{-1}a_6c_6, f_3, g_7=a_5^{-1}, g_8=a_6^{-1}b_7^{-1}$ and $g_9= a_8c_6c_7 ^{-1}, f_4=b_8, g_{10}=a_8^{-1}.$
Then the representation of 
$b_{1\to 6}b_8$  in terms of the generators of $M$ has the form
\begin{equation}
\label{eq:11}b_{1\to 6}b_8=g_1^{\zeta}g_2^{\zeta '}f_1g_3^{\zeta}g_4^{\zeta '}f_2g_5^{\zeta}g_6^{\zeta '}f_3g_7^{\zeta}g_8^{\zeta '}g_9^{\zeta ''}f_4g_{10}^{\zeta ''}.\end{equation} The last three generators commute with each of the first $10$ generators. For given positive $\zeta , \zeta ', \zeta ''$ the equality  
$\zeta \cdot \zeta ' = \zeta ''$ is necessary and sufficient for the indicated occurrence of the element $b_{1\to 6}b_8$ in the submonoid $M$. For null value of $\zeta , \zeta '$ or $\zeta ''$, you can also assume that the corresponding generator  is located as indicated.
\end{lemma}
Proof.
It is clear that the configuration of the factors in (\ref{eq:10}) is preserved for their counterparts in (\ref{eq:11}).
Thus, the $7$th component of the product of these analogs in (\ref{eq:11}) is equal to $c_1^{\zeta_1\zeta_2}$. This value must cancel out due to the $c_1^{-\zeta_3}$ multiplier present in the $g_9$ exponent in the representation (\ref{eq:11}).
\hfill $\square$

Let $\mathbb{H}^8 = \mathbb{H}(1) \times \ldots \times \mathbb{H}(8).$ The scheme of the exact location of the components of the generators of the submonoid $M$ when expressing the element $b_{1\to 6}b_8$ is as follows (empty positions correspond to trivial elements).

\noindent 
 $
 \mathbb{H}^8: b_{1\to 6}b_8=$
 
 $$
 \begin{vmatrix} b_{1\to 6}b_8  =  g_1^{\zeta}g_2^{\zeta '}&f_1&g_3^{\zeta}g_4^{\zeta '}&f_2&g_5^{\zeta}g_6^{\zeta '}&f_3&g_7^{\zeta}g_8^{\zeta '}&g_9^{\zeta ''}&f_4&g_{10}^{\zeta ''}\\
\mathbb{H}(1):  b_1\hspace{0.4cm} = a_1^{\zeta}c_1^{\zeta '}&b_1&a_1^{-\zeta}&&&&&&&\\
\mathbb{H}(2): b_2\hspace{0.4cm}=a_2^{\zeta '}c_2^{\zeta '}&b_2&a_2^{-\zeta '}&&&&&&&\\
\mathbb{H}(3):b_3=&&a_3^{\zeta}c_3^{\zeta}&b_3&a_3^{-\zeta}&&&&&\\
\mathbb{H}(4):b_4=&&a_4^{\zeta '}c_4^{\zeta '}&b_4&a_4^{-\zeta '}&&&&&\\
\mathbb{H}(5):b_5=&&&&a_5^{\zeta}c_5^{\zeta}&b_5&a_5^{-\zeta}&&&\\
\mathbb{H}(6):b_6=&&&&a_6^{\zeta '}c_6^{\zeta '}&b_6&a_6^{-\zeta '}&&&\\
\, \, \, \, \, \, \, \, \, \mathbb{H}(7):1\; \, \, = a_7^{\zeta}&&b_7^{\zeta '}&&a_7^{-\zeta}&&b_7^{-\zeta '}&c_7^{\zeta ''}&&\\
\mathbb{H}(8): b_8= &&&&&&&a_8^{\zeta ''}c_8^{\zeta ''}&b_8&a_8^{-\zeta ''}\\
 \end{vmatrix}
 $$
  \begin{center} Table 4.
 \end{center}

\section{Choosing a direct power of the Heisenberg group and constructing a submonoid in it for which the membership problem is equivalent to the solvability of the given Diophantine equation}
\label{sec:4}

First, a Diophantine  equation (\ref{eq:1}) is taken. Then the equivalent nonnegative Skolem system $S(\upsilon )$ is constructed from this equation. The variables and equations of this system are ordered and written as specified in (\ref{eq:3}--\ref{eq:6}).

To the resulting system $S_{\upsilon}$ we associate the group $\widetilde{\mathbb{H}}=\mathbb{H}^{8e+4d+q+1}$.  We construct a submonoid $M=M$ of the group $\widetilde{\mathbb{H}}$ by defining its generating elements $g_i$ for $i = 1, \ldots , 10e, 10e+1, \ldots , 10e + 6d$  and  $f_j$ for $j = 1, \ldots , 4e, 4e+1, \ldots , 4e+d$ in accordance with the lemmas \ref{le:3} and \ref{le:5}. The form of these generators are defined below. 

\paragraph{\bf 
Construction of submonoid generators associated with the system (\ref{eq:3}).}
 Consider first the equations of system (\ref{eq:3}). For each of these $e$ equations, we sequentially define the corresponding block $\mathbb{H}(i) =\mathbb{H}^8, i = 1, \ldots , e$. 
 Consistently compose the group  
 \begin{equation}
 \label{eq:12}
 \mathbb{H}^{8e} = \prod_{i=1}^e\mathbb{H}(i)
 \end{equation} 
 \noindent  from the obtained blocks. We assume that the group $\mathbb{H}^{8e}$  consists of the first $8e$ factors of the group $\widetilde{\mathbb{H}}$. Let $\bar{M}_1$  be the projection of $M$ into $\mathbb{H}^{8e}$. 
 We define successively $14e$ projections $\bar{g}_i\, (i=1, \ldots ,10e), \bar{f}_j, \, (j = 1, \ldots , 4e)$ of the generators $g_i (i=1, \ldots ,10e), f_j (j = 1, \ldots , 4e)$ of the submonoid $M$ into  the group $\mathbb{H}^{8e}$ with $14$ generators for each block ($10$ generators $\bar{g}_i$ and $4$ generators $\bar{f}_j$).  All these projections generate $\bar{M}_1$, the projection of $M$ into $\mathbb{H}^{8e}.$ The remaining generators of $M$ have trivial projections. 
 
 The generators of the block $\mathbb{H}(1)$ are constructed  in exactly the same way as in Lemma \ref{le:5}. 
 They are $\bar{g}_1, \ldots \bar{g}_{10}$ and $\bar{f}_1, \ldots , \bar{f}_4$. They are considered as elements of the group $\mathbb{H}^{8e}$. All their other  components in $\mathbb{H}^{8e}$ are trivial. The generating elements of the remaining components are determined in the same way. 
 The generators of the submonoid $\bar{M}_1$ in the block  $\mathbb{H}(i)$ (for simplicity we denote $t_i = 10(i-1)$ and $s_i=4(i-1)$) are:
 $$
  \bar{g}_{1+t_i}=a_{1+t_i}c_{1+t_i}a_{7+t_i}, \bar{g}_{2+t_i} =a_{2+t_i}c_{2+t_i}, \bar{g}_{3+t_i}= a_{1+t_i}^{-1}a_{3+t_i}c_{3+t_i},$$ $$\bar{g}_{4+t_i}= a_{2+t_i}^{-1}a_{4+t_i}c_{4+t_i}b_{7+t_i}, \bar{f}_{1+s_i}= b_{1+s_i}b_{2+s_i},$$ 
  \begin{equation}
 \label{eq:13}\bar{f}_{2+s_i}= b_{3+s_i}b_{4+s_i}, \bar{g}_{5+t_i}= a_{3+t_i}^{-1}a_{5+t_i}c_{5+t_i}a_{7+t_i}^{-1},
   \bar{g}_{6+t_i}= a_{4+t_i}^{-1}a_{6+t_i}c_{6+t_i},\end{equation}
   $$ \bar{f}_{3+s_i}=b_{5+s_i}b_{6+s_i}, \bar{g}_{7+t_i}= a_{5+t_i}^{-1}, \bar{g}_{8 +t_i}= a_{6+t_i}^{-1}b_{7+t_i}^{-1},$$ 
  $$\bar{g}_9=a_{8+t_i}c_{6+t_i}c_{7+t_i}^{-1}, \bar{f}_{4+s_i}= b_{8+s_i}, \bar{g}_{10+t_i}=a_{8+t_i}^{-1}.$$

   The analogue of the formula (\ref{eq:11}) for the block 
 $\mathbb{H}(i)$ is the following  formula:
$$
 b_{(1+8(i-1))\to (6+8(i-1))}b_{8i}=\bar{g}_{1+t_i}^{\zeta_{1+t_i}}\bar{g}_{2+t_i}^{\zeta_{2+t_i}}\cdot $$
 \begin{equation}
 \label{eq:14}
 \bar{f}_{1+s_i}
 \bar{g}_{3+t_i}^{\zeta_{1+t_i}}\bar{g}_{4+t_i}^{\zeta_{2+t_i}}\bar{f}_{2+s_i}\bar{g}_{5+t_i}^{\zeta_{1+t_i}}\bar{g}_{6+t_i}^{\zeta_{2+t_i}}\cdot \end{equation}
 $$\bar{f}_{3+s_i}\bar{g}_{7+t_i}^{\zeta_{1+t_i}}\bar{g}_{8+t_i}^{\zeta_{2+t_i}}\bar{g}_{9+t_i}^{\zeta_{3+t_i}}\bar{f}_{4+s_i}\bar{g}_{10+t_i}^{\zeta_{3+t_i}}.
 $$
 
  Just as in Lemma \ref{le:5}, we can conclude that the element $b(1)=b_{1\to 6}b_8 \cdot b_{9\to14}b_{16} \cdot \ldots \cdot b_{(8e-7) \to (8e-2)}b_{8e}$  belongs to $\bar{M}_1$ if and only if all the relations of system (\ref{eq:4}) are satisfied.

\bigskip 
\paragraph{\bf 
Construction of submonoid generators associated with the system (\ref{eq:4}).}
 Consider  the equations of system (\ref{eq:4}). For each of these $d$ equations, we sequentially define the corresponding block $\mathbb{H}(j) =\mathbb{H}^4, j = e+1, \ldots , e+d$. 
 Consistently compose the group  
 \begin{equation}
 \label{eq:15}
 \mathbb{H}^{4d} = \prod_{j=1}^{d}\mathbb{H}(e+j)
 \end{equation} \noindent  from the obtained blocks.  We assume that the group $\mathbb{H}^{4d}$  consists of the  $4d$ factors of the group $\widetilde{\mathbb{H}}$ following after the previously considered factors of $\mathbb{H}^{8e}$ of the group $\widetilde{\mathbb{H}}$.  Now we consider the group $\mathbb{H}^{8e}\times \mathbb{H}^{4d}$, which consists of the first $8e+4d$ factors $\mathbb{H}$ of the group $\widetilde{\mathbb{H}}$.
 
 Let $\bar{M}_2$ be the projection of $M$ into $\mathbb{H}^{4d}$. The projections of the generators of the submonoid $M$ considered above into the components $\mathbb{H}^{4d}$  are trivial. We define successively $7d$ projections $\bar{g}_i, \bar{f}_j$ of the generators $g_i, f_j$ for $i=10e+1, \ldots , 10e+6d$ and  $j = 4e+1, \ldots , 4e+d$ of $M$   ($6$ projections $\bar{g}_i$ and $1$ projection of $\bar{f}_j$ for each block $\mathbb{H}^4$)   into $\mathbb{H}^{4d}$ which are the generators of the submonoid $\bar{M}_2$ of $\mathbb{H}^{4d}$.  Recall that the projections of these generators into $\mathbb{H}^{8e}$ are trivial. Therefore, the projection $\bar{M}_{1,2}$ of the submonoid $M$ into the group $\mathbb{H}^{8e+4d}$ coincides with $\bar{M}_1\bar{M}_2 $.
 
 The generators of the block $\mathbb{H}(e+j), j = 1, \ldots , d,$ are constructed  in exactly the same way as in Lemma \ref{le:3}. The generators of the submonoid $\bar{M}_2$ in the block  $\mathbb{H}(e+j)$ are (for simplicity we denote $r_j = 4(j-1)$): 
 
$$ \bar{g}_{10e+1+r_j}=a_{10e+1+ r_j}c_{10e + 1 +r_j}c_{10e+4+r_j}, \bar{g}_{10e+2+r_j}=a_{10e+2+ r_j}c_{10e+2+r_j}c_{10e+4+r_j},$$
 \begin{equation}
 \label{eq:16} \bar{g}_{10e+3+r_j}=a_{10e+3+r_j}c_{10e+3+r_j}c_{10e+4+r_j}^{-1}, \bar{g}_{10e+4+r_j}=a_{10e+1+r_j}^{-1}, \end{equation}
  $$ \bar{g}_{10e+5+r_j}=a_{10e+2+r_j}^{-1}, \bar{g}_{10e+6+r_j}= a_{10e+3+r_j}^{-1}, \bar{f}_{4e+j}=b_{4e+1+j}b_{4e+2+j}b_{4e+ 3+j}$$ 
 \noindent with trivial other components. 
 
 The analogue of the formula (\ref{eq:9}) for the block 
 $\mathbb{H}(i+j)$ is the following  formula:
 \begin{equation}
 \label{eq:17}
 b_{e+1+r_j\to e+3+r_j}=\bar{g}_{e+1+r_j}^{\zeta_{e+1+r_j}}\bar{g}_{e+2+r_j}^{\zeta_{e+2+r_j}}\bar{g}_{e+3+r_j}^{\zeta_{e+3+r_j}}\bar{f}_{e+1+r_j}\bar{g}_{e+4+r_j}^{\zeta_{e+1+r_j}}\bar{g}_{e+5+r_j}^{\zeta_{e+2+r_j}}\bar{g}_{e+6+r_j}^{\zeta_{e+3+r_j}}.
\end{equation} 
 
 Just as in Lemma \ref{le:3}, we can conclude that the element $b(2)=b_{(10e+1)\to (10e+3)} \cdot b_{(10e+5) \to (10e+7)} \cdot \ldots \cdot b_{10e+4(d-1) +1 \to 10e+4d-1)}$  belongs to $\bar{M}_2$ if and only if all the relations of system (\ref{eq:4}) are satisfied. 

Thus, we have defined the projections $\bar{g}_i, \bar{f}_j$
 of all $14e+7d$ generating elements $g_i, f_j$ of the submonoid $M$ into the product $\mathbb{H}^{8e}\times \mathbb{H}^{4d}$ of the first $8e+4d$ factors $\mathbb{H}$ of the group $\widetilde{\mathbb{H}}.$
 
\bigskip
Summarizing the above, we conclude that the element $b(1)b(2)$ belongs to the submonoid $\bar{M}_{1,2}$ if and only if the joint system of equations (\ref{eq:3} -\ref{eq:4}) is solvable.  

\bigskip
\paragraph{\bf 
Construction of submonoid generators associated with the systems (\ref{eq:5}) and (\ref{eq:6}).}
Since all the variables of the systems (\ref{eq:3}) and (\ref{eq:5}) are different, both these systems are decidable together. It remains to take into account the equalities between these variables.

The difference between this construction and the above constructions related to systems (\ref{eq:3})  and (\ref{eq:4})  is that the considered direct product of the blocks of the group $\mathbb{H}$  is not only expanded by new factors $\mathbb{H}$, but also the already defined projections of the generators of the submonoid $M$ are supplemented with new components. In other words, these generators are modified by the added components.  

Let us add to the constructed group $\mathbb{H}^{8e+4d}$ by $q+1$ factors $\bar{H}$ and get the group $\widetilde{\mathbb{H}}$, where $q$ is the number of  equations in the system (\ref{eq:5}). Let's assign to these components the numbers $8e+4d+i$ for $i = 1, \ldots , q$ and $8e+4d+q+1$ relatively. 

Then for any equation $P_{e+d+k}\, (k=1, \ldots , q)$ of the form  $\zeta_{i(k)}  = \zeta_{j(k)}$ from (\ref{eq:5}) we add some elements  to the $8e+4d+k$th component of $\widetilde{\mathbb{H}}$ as follows. 

First, for each $k=1, \ldots , q$ we find among the representations (\ref{eq:14}) and (\ref{eq:17}) one of the generating elements $g$ of the $M$ submodule whose projection exponent is equal to $\zeta_{i(k)}$. Add the element $c_{8e+4d+k}$ to the component $8e+4d+ k$ of $g$. Then we will perform a similar operation corresponding to the exponent $\zeta_{j(k)}$. This component will be trivial in the considered product of generating elements of the submonoid $M$ if and only if $\zeta_{i(k)}=\zeta_{j(k)}$. 
 
Then for equation (\ref{eq:6}), we find among the representations (\ref{eq:14}) and (\ref{eq:17}) one of the generating elements $g$ of the $M$ submodule whose projection exponent is equal to $\zeta_{t}$ and add the element $c_{8e+4d+q+1}$ to the $8e+4d+q+1$th component of $g$. Note, that this component is equal to $c_{8e+4d+q+1}^{\zeta_t}$ in the considered product of the generators of $M$, i.e., this product is equal to $b(1)b(2)c_{8e+4d+q+1}^{|\upsilon |}.$  

The process of constructing the generators of the submonoid $M$ of the group $\widetilde{\mathbb{H}}$  is completed.

\section{Main results}
\label{sec:5}
 In this section we   give formal proofs of the main results.

\begin{theorem}
\label{th:1}
For any Diophantine equation (\ref{eq:1}), there exists a direct power $\widetilde{\mathbb{H}}=\mathbb{H}^n$ of the Heisenberg group $\mathbb{H}$, a finitely generated submonoid $M $ in the group $\widetilde{\mathbb{H}}$ and an element $g(\upsilon )\in \widetilde{\mathbb{H}}$ such that the equation (\ref{eq:1}) is solvable in integers if and only if $g (\upsilon )$ belongs to $M.$ The exponent $n$, the element $g(\upsilon )$, and the finite set of generators of the submonoid $M$ are effectively determined. The submonoid $M$ depends only on the Diophantine polynomial $D$ on the left side (\ref{eq:1}).
\end{theorem}
Proof.
Suppose that the equation (\ref{eq:1}) has a solution in integers. From the equation (\ref{eq:1}), we construct a nonnegative Skolem system $S_{\upsilon}$ equivalent to it, as explained in the point 1.1.1 and Lemma \ref{le:1}. 
Suppose that $\zeta_1, \ldots , \zeta_{3e+3d}$ is an integer solution to the system  $S_{\upsilon}$. In this case equations (\ref{eq:3}--\ref{eq:6}) turn into equalities.

Define $n = 8e+4d+q+1$ and the group $\widetilde{\mathbb{H}}=\mathbb{H}^n.$ 
Construct the generating elements $g_i$ for $i = 1, \ldots , 10e, 10e+1, \ldots , 10e + 6d$  and  $f_j$ for $j = 1, \ldots , 4e, 4e+1, \ldots , 4e+d$ of the submonoid $M$ of the group $\widetilde{\mathbb{H}}$, as described in the section \ref{sec:4}. Define the element 
\begin{equation} 
\label{eq:18}
g(\upsilon ) =b(1)b(2)c_{8e+4d+q+1}^{|\upsilon |},
\end{equation}
 \noindent
where $b(1) = b_{1\to 6}b_8 \cdot b_{9\to14}b_{16} \cdot \ldots \cdot b_{(8e-7) \to (8e-2)}b_{ 8e}, b(2)=b_{(10e+1)\to (10e+3)} \cdot b_{(10e+5) \to (10e+7)} \cdot \ldots \cdot b_{10e+4(d-1) +1 \to 10e+4d-1)}$. 

  The section \ref{sec:4} shows that the element $b(1)b(2)$ is represented in a certain way as the product of the projections of the generating elements of the submonoid $M$ onto the group $\mathbb{H}^{8e+4d} $ (formulas (\ref{eq:14} and (\ref{eq:17})). Consider the corresponding product of generating elements of the submonoid $M$. It follows from their construction and the fulfillment of the equalities (\ref{eq:5}) that all components $8e+4d+k\, (k=1, \ldots , q)$ of this product are trivial. The $(8e+4d+q+1)$th component in view of the equality (\ref{eq:6}) is equal to $c_{8e+4d+q+1}^{|\upsilon |}$. Hence, the element $g(\upsilon )$ belongs to the submonoid $M$.

Suppose now that the element $g(\upsilon )$ belongs to the submonoid $M$. Its projection $b(1)b(2)$ is represented in a certain way as a product of the projections of generating elements of the submonoid $M$ only if the equalities (\ref{eq:3}) and (\ref{eq:4}) hold ( Lemmas \ref{le:3} and \ref{le:5}, formulas (\ref{eq:14}) and (\ref{eq:17})). This product completely determines the product of generators of the submonoid $M$. In this case, the components with numbers $8e+4d+k\, (k=1, \ldots , q)$ must be trivial, which corresponds to the fulfillment of the equalities (\ref{eq:5}). The ($8e+4d+q+1$)th component must be equal to $c_{8e+4d+q+1}^{|\upsilon |},$ which means that (\ref{eq:6}) satisfied. Consequently, the exponents $\zeta_1, \ldots , \zeta_{3(e+d)}$, with which the generators of the submonoid $M$ enter the representation of the element $g(\upsilon )$, are the solution of the system $S_{\upsilon }$.
\hfill $\square$

Recall that Hilbert's 10th problem is the question of the existence of an algorithm that, given a
Diophantine equation determines whether it has an integer solution.
Yu.V. Matiyasevich (see \cite{Mat1}--\cite{MatRob}) proved that such an algorithm does not exist. In addition, he established that there exists a Diophantine polynomial $D_0(\zeta_1, \ldots , \zeta_{t})$ with a zero constant term such that there is no algorithm that determines the solvability of equations of the form
\begin{equation}
\label{eq:19}
D_0(\zeta_1, \ldots , \zeta_t) = \upsilon ,\, \upsilon \in \mathbb{Z}.
\end{equation}

From the undecidability of Hilbert's 10th problem and  Theorem \ref{th:1}, it follows that the submonoid membership problem in the class of finite direct powers of the Heisenberg group is undecidable.

The existence of an algorithmically unsolvable equation of the form (\ref{eq:19}) with a fixed left-hand side and parameter $\upsilon$ allows us to establish the following stronger assertion.

\begin{theorem}
\label{th:2}
For sufficiently large $n\in \mathbb{N}$, the direct power $\widetilde{\mathbb{H}}=\mathbb{H}^n$ of the Heisenberg group $\mathbb{H}$ contains a finitely generated submonoid $M$
 with an unsolvable membership problem.
\end{theorem}
Proof. First, an equation of the form (\ref{eq:19}), which is unsolvable in integers, is taken. Then the equivalent nonnegative Skolem system $S(\upsilon )$ is constructed from this equation. The rest of the proof completely repeats the proof of the Theorem \ref{th:1}.  Variations of the parameter $\upsilon$ in the equation (\ref{eq:19}) correspond to variations of the element $g(\upsilon ).$ The submonoid $M$ does not change. An element $g(\upsilon )$ belongs to the submonoid $M$ if and only if the system $S_{\upsilon}$ is solvable in 
nonnegative integers. This is equivalent to saying that the equation (\ref{eq:19}) with this parameter is solvable in integers. This implies the assertion of the theorem.
\hfill $\square$

Note that the existence of a finitely generated submonoid with an unsolvable occurrence problem in a finitely generated nilpotent group implies the existence of a similar submonoid in the corresponding free nilpotent group.

\begin{proposition}
\label{pr:1}
For $k, c \in \mathbb{N}$, let $N$ be a  $k$-generated nilpotent group of class $c$ that has a finitely generated submonoid $M$ with an undecidable membership problem. Then the free nilpotent group $N_{k,c}$ of rank $k$ of the class $c$ contains a finitely generated submonoid $\widetilde{M}$ with an undecidable membership problem.
\end{proposition}
Proof. Consider the natural homomorphism $\mu : N_{k,c} \rightarrow N.$ Let $\widetilde{M}$ denote the full pre-image of the submonoid $M$ in $N_{k,c}$. An element $g\in N$ belongs to $M$ if and only if any of its inverse images $\widetilde{g}$ belongs to $\widetilde{M}$. It remains to note that the submonoid $\widetilde{M}$ is finitely generated.

Let $M$ is generated by elements $g_1, \ldots , g_l$. For each of these generators $g_i$, take some inverse image $\tilde{g}_i$ in the group $N_{k,c}.$ The group $N_{k,c}$ is Noetherian, so ker($\mu$) is a finitely generated subgroup. Let ker($\mu$) = gp($f_1, \ldots , f_t$). Then the submonoid $\widetilde{M}$ is generated by the elements $\tilde{g}_1, \ldots , \tilde{g}_l, f_1^{\pm 1}, \ldots , f_t^{\pm 1}$. 
\hfill $\square$

It follows from the Theorem \ref{th:2} and the Proposition \ref{pr:1} that a submonoid with an unsolvable membership problem exists in any free nilpotent group $N_{k,c}$ for $c\geq 2$ of sufficiently large rank $ k$.

\end{document}